\documentclass[11pt,reqno]{amsart}
\usepackage[dvipdfm]{graphicx, color}
\usepackage{amssymb, amsmath,amsthm}
\newtheorem{theorem}{{\bf Theorem}}
\newtheorem{lemma}[theorem]{{\bf Lemma}}
\newtheorem{corollary}[theorem]{{\bf Corollary}}

\newtheorem{reduction}[theorem]{Reduction Hypothesis}
\numberwithin{equation}{section}

\begin{document}

%
%

\title[Entropy, Weil-Petersson translation distance and 
Gromov norm]{Entropy, Weil-Petersson translation distance and \\ 
Gromov norm for surface automorphisms}

\author[S. Kojima]{%
    Sadayoshi Kojima
}
\address{%
        Department of Mathematical and Computing Sciences \\
        Tokyo Institute of Technology \\
        Ohokayama, Meguro \\
        Tokyo 152-8552 Japan
}
\email{%
        sadayosi@is.titech.ac.jp
}



\thanks{%
The author is partially supported by Grant-in-Aid for Scientific Research (A)
 (No. 22244004), JSPS, Japan
} 

\begin{abstract} 
Thanks to a theorem of Brock on comparison of 
Weil-Petersson translation distances and hyperbolic volumes of mapping tori for pseudo-Anosovs, 
we prove that the entropy of a surface automorphism in general has linear 
bounds in terms of Gromov norm of its  mapping torus 
from below and in bounded geometry case from above.  
We also prove that the Weil-Petersson translation distance does the same from both sides in general. 
The proofs are in fact immediately derived from the theorem of Brock together with 
some other strong theorems and small observations.   
\end{abstract}

\maketitle

%
%

\section{Introduction}\label{Sect:introduction}

Let  $\varSigma = \varSigma_{g,n}$ be an orientable surface of genus $g$ with $n$ punctures, 
$\mathcal{T} = \mathcal{T}_{g,n}$  the Teichm\"uller space of  $\varSigma$    
and  $\mathcal{M} = \mathcal{M}_{g,n}$  the (orientation preserving) 
mapping class group of  $\varSigma$.   
An element of  $\mathcal{T}$  is defined by  
a pair  $(R, f)$  of a Riemann surface  $R$  homeomorphic to  $\varSigma$  
with a marking homeomorphism  $f : \varSigma \to R$.   
A mapping class  $\varphi \in \mathcal{M}$   naturally acts 
on  $(R, f) \in \mathcal{T}$ by  $(R, f \circ \varphi^{-1})$.  
Throughout this paper, 
we assume that  $3g - 3 + n \geq 1$  or  $(g, n) = (0,3)$, 
and fix  $g$  and  $n$.  
Also, 
we omit a marking homeomorphism  $f$  to indicate 
a marked Riemann surface in  $\mathcal{T}$.  

In this paper, 
we discuss three invariants measuring the complexity of  
a conjugacy class of a mapping class  $\varphi \in \mathcal{M}$.  
One is the infimum of topological entropies of 
surface automorphisms in  $\varphi$,   
which we denote by  ${\rm ent} \, \varphi$.    
As we will point out in Corollary \ref{Cor:Entropy=T}, 
it is equal to the Teichm\"uller translation distance 
of  $\varphi \in \mathcal{M}$, 
\begin{equation*}
	|| \varphi ||_{\rm T} = \inf_{R \in \mathcal{T}} d_{\rm T}(R, \varphi(R)),     
\end{equation*}
where  $d_{\rm T}$  denotes the Teichm\"uller distance on  $\mathcal{T}$.  
The second one is the Weil-Petersson translation distance defined by 
\begin{equation*}
	|| \varphi ||_{\rm WP} = \inf_{R \in \mathcal{T}} d_{\rm WP}(R, \varphi(R)),     
\end{equation*}
where  $d_{\rm WP}$  denotes the Weil-Petersson distance on  $\mathcal{T}$.  
The last one through three dimensional manifolds is 
the Gromov norm   $|| \cdot ||_{{\rm Gr}}$  of 
the mapping torus of  $\varphi$, 
\begin{equation*} 
	N_{\varphi} =\varSigma \times [0,1]/ (x, 1) \sim (\varphi(x), 0).  
\end{equation*} 
Here the notation  $N_{\varphi}$  indicates only a mapping class because 
the topological type of  $N_{\varphi}$  
does not depend on the choice of representatives of  $\varphi$.  
In fact, 
it depends only on the conjugacy class of  $\varphi$.  
Hence any topological invariant of  $N_{\varphi}$  such as 
the Gromov norm is an 
invariant of a mapping class  $\varphi$.  

The starting point of this paper is a comparison of 
the Weil-Petersson translation distances of pseudo-Anosovs and 
volumes of their mapping tori by Brock in \cite{Brock}.  
We discussed similar inequalities in Theorem 3.2 of  \cite{KKT}, 
where we replaced the Weil-Petersson translation distance in 
Brock's work by the entropy.  
Here we generalize it further to one for general surface automorphisms  
as follows. 

\begin{theorem}\label{Thm:1} 
There is a constant  $A = A(g, n) > 0$  which depends only on 
the genus  $g$  and the number  $n$  of punctures  
such that the inequality 
\begin{equation*}
	A \, || N_{\varphi} ||_{{\rm Gr}} \leq {\rm ent} \, \varphi 
\end{equation*} 
holds for any  $\varphi \in \mathcal{M}_{g,n}$.  
Furthermore for any  $\varepsilon > 0$,  
there is a constant  $B = B(g, n, \varepsilon) > 0$  which 
depends only on  $g, \, n$  and  $\varepsilon$  such that 
\begin{equation*} 
	{\rm ent} \, \varphi \leq B \, || N_{\varphi} ||_{{\rm Gr}}  
\end{equation*} 
holds for any  $\varphi \in \mathcal{M}_{g,n}$  of which the hyperbolic pieces of 
$N_{\varphi}$  contain no closed geodesics of length  $< \varepsilon$.  
\end{theorem}  

We also generalize the theorem of Brock to one for 
general surface automorphisms as follows.  

\begin{theorem}\label{Thm:2}
There exists a constant  $C = C(g, n)$  which depends only on  $g$  and  $n$  
such that 
\begin{equation*}
	C^{-1} \, || \, N_{\varphi} ||_{{\rm Gr}} 
	\leq || \varphi ||_{{\rm WP}} 
	\leq C \, || \, N_{\varphi} ||_{{\rm Gr}} 
\end{equation*} 
holds for any  $\varphi \in \mathcal{M}_{g,n}$.   
\end{theorem}

Both theorems are immediately derived from Brock's theorem together with 
some existing strong theorems and small observations.  
We collect necessary materials first with reasonable account 
of explanations as packages in the next section,  
and then present reducing argument in the last section.  
For the reader's convenience, 
we duplicate some of the arguments in  \cite{KKT}.  
\medskip 

\noindent 
{\bf Acknowledgment :} 
I would like to thank Eiko Kin and Mitsuhiko Takasawa 
for their helpful comments.

%
%

\section{Packages}


\subsection{Entropy} 

Let us review the topological entropy defined by 
Adler, Konheim and McAndrew in  \cite{AKM}  
for a self map  $f : X \to X$  of a compact topological space  $X$.  
Choose an open covering  $\mathcal{U}$  of  $X$,  
and 
let  $N(m, \mathcal{U})$  be the minimal number of members 
in a refined covering  
$\mathcal{U} \vee f^{-1}(\mathcal{U}) \vee \cdots \vee f^{-(m-1)}(\mathcal{U})$  
by intersections  
that cover  $\varSigma$.  
Then  
\begin{equation*}
	h(f) = \sup_{\mathcal{U}} \left( \limsup_{m \to \infty} \frac{1}{m} \log N(m, \mathcal{U}) \right)
\end{equation*}
is called the topological entropy of  $f$.  
It is shown in  \cite{AKM}  to always exist and 
to satisfy the identity  
\begin{equation*} 
	h(f^n) = n \, h(f)
\end{equation*}   
for any positive integer  $n$.  

Our topological space  $\varSigma$  may not be compact, 
however if we suppose that the covering  $\mathcal{U}$  contains  
open sets including each puncture in their interior, 
then it can be handled as a compact space and 
every argument in  \cite{AKM}  works equally well.  
For a mapping class  $\varphi \in \mathcal{M}$, 
we define the entropy of  $\varphi$  by 
\begin{equation*}
	{\rm ent} \, \varphi = \inf_{f \in \varphi} h(f).  
\end{equation*}
In particular, 
we have the identity 
\begin{equation}\label{Eq:Entropy} 
	{\rm ent} \, \varphi^n = n \, {\rm ent} \, \varphi. 
\end{equation}


\subsection{Thurston's classification of surface automorphisms}  

According to a seminal work of Thurston in \cite{Thurston2} (see also \cite{FLP}), 
the elements of  $\mathcal{M}$  are classified into three types : 
periodic, pseudo-Anosov and reducible.  
This classification has become standard, 
but it is better to note that the classification is not in trichotomy.  
In fact, 
there is a class which is periodic and also reducible.  

A pseudo-Anosov mapping class  $\varphi$  is known to have a good deal of structures.  
It has a nice representative  $f : \varSigma \to \varSigma$  which 
leaves a pair of stable and unstable foliations invariant.  
The expanding factor  $\lambda$  of the unstable foliation is called 
the dilatation of  $f$. 
Thurston \cite{Thurston2} (see also \cite{FLP})  showed that 
if  $\varphi \in \mathcal{M}$  is pseudo-Anosov, 
then 
\begin{equation*}
	{\rm ent} \, \varphi = \log \lambda.  
\end{equation*}


\subsection{Bers' classification of surface automorphisms} 

Bers gave an alternative proof of Thurston's classification in  \cite{Bers}   
by adopting extremal approach in quasiconformal mappings. 
Since Bers' classification is more suitable for us, 
we briefly review it. 

Bers calls  $\varphi$  elliptic if  $|| \varphi ||_{\rm T} = 0$  and 
there is a Riemann surface   $R \in \mathcal{T}$  such that  $\varphi(R) = R$.  
This is identical with Thurston's periodic case.    
Call  $\varphi$  hyperbolic if  $|| \varphi ||_{\rm T} >0$  and there is 
a Riemann surface  $R \in \mathcal{T}$  which attains the infimum.  
This is also identical with Thurston's pseudo-Anosov case.    
The remaining cases in Bers' classification are called 
parabolic if  $|| \varphi ||_{\rm T} = 0$  or pseudo-hyperbolic if  $|| \varphi ||_{\rm T} > 0$,  
where there are no Riemann surfaces in  $\mathcal{T}$  which attain the infimum for 
either cases.  
These are reducible but not periodic in Thurston's classification. 

Then Bers showed that if  $\varphi$  is hyperbolic (= pseudo-Anosov), 
then there is a unique invariant Teichm\"uller geodesic  $\ell_{\rm T}(\varphi)$  
for the action of  $\varphi$  on  $\mathcal{T}$,  
and the identity 
\begin{equation*}
	\log \lambda = d_{\rm T}(R, \varphi(R)) = || \varphi ||_{\rm T} 
\end{equation*} 
holds for any  $R \in \ell_{\rm T}(\varphi)$.


\subsection{Gromov norm} 

According to Section 6.5 of \cite{Thurston1}, 
we recall the definition of the Gromov norm for a 3-manifold with toral boundary.  
Let   $M$  be a compact orientable 3-manifold whose boundary consists of tori.  
Then, 
the Gromov norm is defined to be the limit infimum of 
$\ell_1$-norm of cycles representing the relative fundamental class  
$[M, \partial M] \in H_3(M, \partial M; \mathbb{R})$  with 
vanishing boundary condition, 
namely, 
\begin{equation*}
	||[M||_{{\rm Gr}} = 
	\liminf_{\varepsilon \to 0} \, \{ ||c||_1 \, ; \, 
		[c] = [M, \partial M], \text{and} \; || \partial c || \leq \varepsilon \}.  
\end{equation*}

In \cite{Soma}, 
Soma proved the splitting theorem for the Gromov norm 
for compact orientable 3-manifolds whose boundary consists of tori.  
His theorem,  
together with the solution of the geometrization conjecture 
by Perelman  \cite{Perelman1, Perelman2}  
establish a simpler statement than the original in  \cite{Soma} as follews.  
Let  $v_3 = 1.01494...$  be the volume of the regular ideal hyperbolic 3-simplex. 

\begin{theorem}[Soma \cite{Soma}]\label{Thm:Soma}
The Gromov norm of a compact orientable 3-manifold  $M$   
whose boundary consists of tori is equal to  $v_3$  times 
the sum of volumes of hyperbolic pieces in 
the geometric decomposition of  $M$. 
\end{theorem}

\noindent
Thus the term  $|| N_{\varphi} ||_{{\rm Gr}}$  in the inequalities of the theorems   
can be replaced by the sum of volumes of hyperbolic parts 
of  $N_{\varphi}$  in the geometric decomposition.  

One of useful properties of the Gromov norm for us is that 
if  $\widetilde{M} \to M$  is an $n$-fold covering, 
then  
\begin{equation}\label{Eq:GromovNorm} 
	|| \widetilde{M} ||_{{\rm Gr}} = n || M ||_{{\rm Gr}}.  
\end{equation}  
This fact will be used to reduce the problem in question in the next subsection.


\subsection{Reduction}  

To each mapping class  $\varphi \in \mathcal{M}$,  
is associated a reduced maximal family  $\mathcal{C}$  of disjoint 
simple closed curves on  $\varSigma$  which is invariant 
by the action of  $\varphi$.  
This means more precisely that, 
$\varSigma - \mathcal{C}$  is divided into two groups  $U$  and  $V$  
which consist of 
components of  $\varSigma - \mathcal{C}$    
so that  $\varphi$ restricts to  $U$  is periodic and  to  $V$  pseudo-Anosov, 
where  $U$  and  $V$  may be disconnected.  

We know from  (\ref{Eq:Entropy})  and  (\ref{Eq:GromovNorm})  that   
${\rm ent}$  and  $|| \; ||_{\rm Gr}$  behave well under taking powers.  
The same identities hold for the Weil-Petersson translation distance and Teichm\"uller translation distance 
by definition, 
namely,  
\begin{equation*} 
	|| \varphi^n ||_{{\rm WP}} = n || \varphi ||_{{\rm WP}}  
	\qquad  \text{and}  \qquad 
	|| \varphi^n ||_{{\rm T}} = n || \varphi ||_{{\rm T}}
\end{equation*} 
hold for any positive integer  $n$.  
Thus,
it is sufficient to show the inequalities in the theorems in the introduction   
for some power of  $\varphi$,  
and we may assume by taking some power of  $\varphi$  that 
$\varphi$  (with the same notation)   
preserves each component of  $\varSigma - \mathcal{C}$.  
Then  $\varSigma$  will be a union of subsurfaces 
$S_1 \cup S_2 \cup \cdots \cup S_k$  split by  $\mathcal{C}$, and 
the restriction of  $\varphi$  to  $S_i$  defines a mapping 
class  $\varphi_i$  on each component  $S_i$.  
If some  $\varphi_i$  is periodic, 
again taking farther power, 
we may assume that it is the identity.  
Thus, 
relabeling suffices 
and letting  $j$  be the number of pseudo-Anosov components, 
we assume in the sequel that 

\begin{reduction}\label{Reduction} 
The restriction  $\varphi_i$  of  $\varphi$  to  $S_i$  in  
$\varSigma  - \mathcal{C} = S_1 \cup S_2 \cup \cdots \cup S_k$  
preserves  $S_i$,    
and is pseudo-Anosov if  $i \leq j \, (\leq k)$  and 
the identity if  $ i > j$.  
\end{reduction} 

\noindent 
Note that  $j$  could range from zero to  $k$.  
Also  $k$  could be  $1$.


\subsection{Geometric decomposition of  $N_{\varphi}$}  
 
Thurston proved in \cite{Thurston3} (cf. \cite{OK})  that 
the mapping torus of a pseudo-Anosov admits 
a hyperbolic structure.  
The hyperbolic  structure in dimension $\geq 3$  is known to be 
unique by Mostow rigidity.   
Under the Reduction Hypothesis \ref{Reduction}, 
the geometric decomposition of  $N_{\varphi}$  is quite simple and 
the hyperbolic pieces consist of  $N_{\varphi_i}$'s for  $i \leq j$.  
Thus by Theorem \ref{Thm:Soma} of  Soma,  
we have 
\begin{equation}\label{Eq:Gromov=Vol}
	|| N_{\varphi} ||_{{\rm Gr}} 
	= \sum_{i=1}^j || N_{\varphi_i} ||_{{\rm Gr}}
	= \frac{1}{v_3} \sum_{i=1}^j {\rm vol} \, N_{\varphi_i}.
\end{equation}


\subsection{Teichm\"uller metric} 

The Teichm\"uller distance beween two Riemann surfaces  
$R, R' \in \mathcal{T}$  is originally 
defined by 
\begin{equation*}
	d_{\rm T}(R, R') = \inf_h \sup_{x \in R} \log K_h(x)
\end{equation*} 
where  $h : R \to R'$  is a quasiconformal map, 
$K_h(x)$  is a dilatation of  $h$  at  $x$.  
To see the infinitesimal form of  $d_{{\rm T}}$,  
we recall a little analysis by quoting a few well known results 
from the text book by Gardiner and Lakic  \cite{GL}.   

Let  $R$  be a Riemann surface.   
Also,  
let  $T^{1,0} R$  and  $T^{0,1}R$  be 
the holomorphic part and the anti-holomorphic 
part of the complex cotangent bundle over  $R$  respectively.  
A Beltrami differential  $\mu$  on  $R$, 
which is a section of a line bundle  
$(T^{0,1}R) \otimes (T^{1,0}R)^*$,   
represents an infinitesimal deformation of complex structure of  $R$  
including holomorphically trivial ones 
caused by diffeomorphisms.  
We let  ${\rm L}_{\infty}(R)$  be the space of 
uniformly bounded Beltrami differentials,  
namely  $\{ \mu \, : \, || \mu ||_{\infty} < \infty \}$.  
This is an infinite dimensional space.  

Another object of interest is a quadratic differential  $q$  
which is a section of a line bundle  $(T^{1,0})^{\otimes 2}$.  
Let  $Q(R)$  be the space of holomorphic quadratic differentials 
with bounded  $L_1$  norm  $|| q || = \int_R |q| < \infty$.  
By Riemann-Roch, 
the dimension of  $Q(R)$  is equal to  $3g-3+n$.  
Then, 
there is a natural pairing defined by 
\begin{equation}\label{Eq:Pairing}
	(\mu, q) = \int_R \mu q.  
\end{equation}  
If we let  
\begin{equation*} 
	N = \{ \mu \, : \, (\mu, q) = 0 \; \text{for all} \: q \in Q(R) \}, 
\end{equation*}  
then  $L_{\infty} (R)/N$  can be identified with the tangent space  $T_R \mathcal{T}$   
of  $\mathcal{T}$  at  $R$.  
A standard argument in functional analysis implies that, 

\begin{lemma}[cf. 3.1 Theorem 2 in \cite{GL}]\label{Lem:NaturalParing}
The natural pairing  {\em (\ref{Eq:Pairing})}  induces 
an isomorphism of  $Q(R)^*$  to  $L_{\infty}(R)/N \cong T_R \mathcal{T}$.  
\end{lemma}

\begin{lemma}[cf. 4.12 Theorem 13 in \cite{GL}]\label{Lem:InfinitesimalTeichmuller} 
	The norm  $|| \mu || = \sup \{ (\mu, q) \, ; \, || q || \leq 1\}$  dual to 
	the  $L_1$ norm on  $Q(R) \cong T_{R}^*\mathcal{T}$  is 
	the infinitesimal form of the Teichm\"uller distance  $d_{\rm T}$.    
\end{lemma}


\subsection{Weil-Petersson metric} 

The Weil-Petersson metric comes from the  $L^2$  inner product 
on  $Q(R) \cong T_R^* \mathcal{T}$  defined by 
\begin{equation*}
	\langle q, q'  \rangle_{{\rm WP}} = \int_X \frac{\bar{q} q'}{\rho^2} 
\end{equation*}
where  $q, q'  \in Q(R)$  and 
$\rho(z) |dz|$  is the hyperbolic metric on  $\varSigma$  
conformally equivalent to the complex structure on  $R$.  
The Weil-Petersson metric is defined to be a Riemannian part of 
the dual Hermitian metric to the above co-metric 
on the cotangent space.  

The Weil-Petersson metric is known to have negatively curved 
sectional curvature by the work of several authors, 
see for example Wolpert \cite{Wolpert1}.   
Also it is known to be incomplete by Wolpert in  \cite{Wolpert0}, 
but geodesically convex also by Wolpert in  \cite{Wolpert2} : 
for each pair of points, 
there exists a unique distance realizing joining curve.  

Daskalopoulis and Wentworth \cite{DW} 
showed that if  $\varphi$  itself is pseudo-Anosov, 
then it admits the unique invariant Weil-Petersson geodesic  $\ell_{{\rm WT}}(\varphi)$  
of the action by  $\varphi$   on  $\mathcal{T}$,  
and the identity 
\begin{equation*}
	d_{\rm WP} (R, \varphi(R)) =  || \varphi ||_{\rm WP} 
\end{equation*} 
holds for any  $R \in \ell_{{\rm WP}}(\varphi)$.


\subsection{Teichm\"uller versus Weil-Petersson}

Lemma \ref{Lem:InfinitesimalTeichmuller} says that 
the $L_1$ norm on  $Q(R)$  is the dual norm to the infinitesimal form 
of the Teichm\"uller distance through the natural paring 
established in lemma \ref{Lem:NaturalParing}.  
Then, 
by the Cauchy-Schwarz inequality, 
we have 
\begin{equation*}
	|| q ||^2 
	= \left( \int_R |q| \right)^2 
	= \left( \int_R \rho \cdot \frac{|q|}{\rho} \right)^2 
	\leq \int_R \rho^2 \cdot \int_R \frac{\bar{q} q}{\rho^2} 
	= {\rm Area} \, \varSigma \cdot \langle q, q \rangle_{{\rm WP}}.  
\end{equation*}
Thus, 
this infinitesimal inequality on the dual space implies the inequality 
of the other direction in two distances,  
\begin{equation*}
	(\sqrt{{\rm Area} \, \varSigma})^{-1} \, d_{{\rm WP}} \leq d_{{\rm T}}, 
\end{equation*} 
which was originally proved by Linch in \cite{Linch}. 
This implies 

\begin{lemma}\label{Lem:WPleqT}
	There is a constant  $D = D(g, n)$  which depends only on  $g, n$  such that the inequality, 
	\begin{equation}\label{Eq:WPleqT}
	|| \varphi ||_{{\rm WP}} \leq  D \, || \varphi ||_{{\rm T}},  	
	\end{equation} 
	holds for any  $\varphi \in \mathcal{M}$.  
\end{lemma}
\begin{proof} 
Choose a Riemann surface  $R$  on the Teichm\"uller  geodesic  
$\ell_{{\rm T}}(\varphi)$  of  $\varphi$  in  $\mathcal{T}$.  
Then, 
we have 
\begin{equation*}
	|| \varphi ||_{{\rm T}} 
	= d_{{\rm T}}(R, \varphi(R)) 
	\geq (\sqrt{{\rm Area} \, \varSigma})^{-1} \, d_{{\rm WP}}(R, \varphi(R)) 
	\geq (\sqrt{{\rm Area} \, \varSigma})^{-1} \, || \varphi ||_{{\rm WP}},  
\end{equation*}
where the last inequality holds because   
$R$  may not be on  $\ell_{{\rm WP}}(\varphi)$.   
\end{proof}

On the other hand, 
since the Weil-Petersson metric is incomplete and 
there is a point at infinity which is of a finite distance from an interior point, 
the inequality of the other direction to   (\ref{Eq:WPleqT})  cannot be 
established in general.  
However, 
there is such a case under some bounded geometry condition.  
For any  $\varepsilon > 0$, 
there is a constant  $\delta > 0$  such that,  
if  $\varphi \in \mathcal{M}$  is pseudo-Anosov and 
if a hyperbolic structure on  $N_{\varphi}$  contains 
no closed geodesics of length $< \varepsilon$, 
then both the Teichm\"uller geodesic  $\ell_{{\rm T}}(\varphi)$  and 
the Weil-Petersson geodesic  $\ell_{{\rm WP}}(\varphi)$  are located 
in the subset  $\mathcal{T}^{\delta}$  of  $\mathcal{T}$  consisting of 
hyperbolic surfaces with no closed geodesics of 
length $< \delta$.  
This fact is proved for the Teichm\"uller metric by Minsky in  \cite{Minsky}  and 
for the Weil-Petersson metric by Brock-Masur-Minsky in  \cite{BMM}. 
This implies 

\begin{lemma}\label{Lem:TleqWP}
For any  $\varepsilon > 0$, 
there is a constant  $E = E(g,n, \varepsilon) > 0$  such that
\begin{equation}\label{Eq:TleqWP}
	E \, || \varphi ||_{{\rm T}} \leq || \varphi ||_{{\rm WP}}, 
\end{equation}
for any pseudo-Anosov $\varphi \in \mathcal{M}$  so that  
$N_{\varphi}$  contains no closed geodesics of length  $< \varepsilon$.  
\end{lemma}\
\begin{proof} 
The subset  $\mathcal{T}^{\delta}$ is invariant by the action of  $\mathcal{M}$, 
and the quotient  $\mathcal{T}^{\delta}/\mathcal{M}$  is 
compact by a theorem of Mumford in \cite{Mumford}.  
Hence  both Teichm\"uller and Wei-Petersson metrics on  $\mathcal{T}^{\delta}$  
are pull backs of metrics on the compact space  $\mathcal{T}^{\delta}/\mathcal{M}$, 
and there is a constant  $E = E(g, n, \varepsilon)$  
which depend only on  $g, n$  and  $\varepsilon$  
such that 
\begin{equation*}
	E \, d_{{\rm T}} \leq d_{{\rm WP}}
\end{equation*} 
on  $\mathcal{T}^{\delta}$.  
Choose a Riemann surface  $R$  on the Weil-Petersson geodesic  
$\ell_{{\rm WP}}(\varphi)$  of  $\varphi$  
on  $\mathcal{T}$.  
Then 
\begin{equation*}
	|| \varphi ||_{{\rm WP}} 
	= d_{{\rm WP}}(R, \varphi(R)) 
	\geq E \, d_{{\rm T}}(R, \varphi(R)) 
	\geq E \, || \varphi ||_{{\rm T}}, 
\end{equation*}
where the last inequality holds because   
$R$  may not be on  $\ell_{{\rm T}}(\varphi)$.   
\end{proof}


\subsection{Weil-Petersson translation distance}  

The augmented structure of the completion  $\overline{\mathcal{T}}$   
with respect to the Weil-Petersson metric 
was extensively studied by Masur in  \cite{Masur}.  
Using this together with Wolpert's convexity, 
Masur and Wolf provided a very clear geodesically embedded picture  
of the frontier of  $\overline{\mathcal{T}}$   in Subsection 1.2 of \cite{MW}.    

One important observation there is that 
the intrinsic Weil-Petersson metric defined on the frontier of  $\mathcal{T}$  is 
identical with the induced metric on the metric completion of  $\mathcal{T}$  
with respect to the Weil-Petersson metric.  
In particular, 
if we let  $\mathcal{O}(\mathcal{C})$  be 
the product $\mathcal{T}(S_1) \times \mathcal{T}(S_2) \times \cdots \times \mathcal{T}(S_k)$   
of Teichm\"uller spaces of  $S_i$'s, 
together with the product metric of 
the Weil-Petersson metircs on each Teichm\"uller space  $\mathcal{T}(S_i)$,  
then  $\mathcal{O}(\mathcal{C})$  is isometrically embedded in  $\overline{\mathcal{T}}$  
as the corresponding frontier of  $\mathcal{T}$. 

Recall that Daskalopoulis and Wentworth \cite{DW} 
established the existence of Weil-Petersson geodesic 
for the action of a pseudo-Anosov  $\varphi$   on  $\mathcal{T}$. 
The following lemma concerns with the other case, 
which is discussed earlier in \cite{Wolpert3}.  
We here provide its quick proof for reader's convenience.   

\begin{lemma} 
The following hold.   
\begin{enumerate}
\item 
	If  $j = 0$,  
	then  $\varphi$  fixes  $\mathcal{O}(\mathcal{C})$  and \,  $|| \varphi ||_{{\rm WP}} = 0$.  
\item 
	If  $k \geq 2$  and  $1 \leq j \leq k$, 
	then 
	there is no invariant geodesic of the action by  $\varphi$  in  $\mathcal{T}$, 
	but there is at least one in  $\mathcal{O}(\mathcal{C}) \subset \overline{\mathcal{T}}$, 
	and  
	\begin{equation}\label{Eq:WPnorm}
		|| \varphi ||_{{\rm WP}} = \sqrt{\sum_{i=1}^j || \varphi_i ||_{{\rm WP}}^2}.  
	\end{equation} 
	The invariant geodesic is not unique if and only if 
	$j < k$  and there exist  $i > j$  so that 
	$S_i$  is not a thrice punctured sphere.  
\end{enumerate}
\end{lemma}

\begin{proof}  
(1)  is clear.  

To see (2), 
choose a point  $R_i  \in \mathcal{T}(S_i)$  on the invariant geodesic of the 
action by  $\varphi_i$  for  $i \leq j$  and any point   $R_i \in \mathcal{T}(S_i)$  if  
$i > j$.  
Then the point  $R = (R_1, R_2, \cdots, R_k)$  on  $\mathcal{O}(\mathcal{C})$  will be 
on the geodesic  $\ell$  in  $\mathcal{O}(\mathcal{C}) \subset \overline{\mathcal{T}}$  
invariant by the action of  $\varphi$.   
The identity (\ref{Eq:WPnorm})  will be obvious by the structure of the metric 
on  $\mathcal{O}(\mathcal{C})$.  

To see non-existence of an invariant geodesic of the action by  $\varphi$  in $\mathcal{T}$,  
suppose to the contrary that we have one  $\ell'$  lying in  $\mathcal{T}$.  
Since  $\mathcal{T}$  with the Weil-Petersson metric has negative sectional 
curvature and is geodesically convex, 
it is a  ${\rm CAT}(0)$ space.  
Then by Chapter II, Corollary 3.11 in \cite{BH}, 
the metric completion  $\overline{\mathcal{T}}$  is also 
a ${\rm CAT}(0)$ space.  
Since  $\ell$  and  $\ell'$  both are invariant by the action of  $\varphi$, 
it is easy to see that they provide asymptotic geodesics, 
namely, 
there are a constant  $K$  and arclength parameterizations for  $\ell$  and  $\ell'$  
so that  $d_{{\rm WP}}(\ell(r), \ell'(t)) < K$  for all  $t \in \mathbb{R}$.     
This contradicts the flat strip theorem, 
Chapter II, Theorem 2.13 in \cite{BH}, 
since  $\ell'$  is located in the part with negative sectional curvature.  
\end{proof}


\subsection{Teichm\"uller translation distance versus Entropy}

By applying Theorem 4 in  \cite{AKM}  to our setting, 
we obtain the identity, 
\begin{equation}\label{Eq:EntropyIsMax}
	{\rm ent} \, \varphi = \max_{1 \leq i \leq j} {\rm ent} \, \varphi_i.  
\end{equation} 
This is the only result in this subsection which we will use in the next section, 
however, 
we would like to add a few more lines that would be worth noting.  

Theorem 9 in \cite{Bers} can be simply stated in our setting by  
\begin{equation*}
	|| \varphi ||_{{\rm T}} = \max_{1 \leq i \leq j} || \varphi_i ||_{{\rm T}}.  
\end{equation*}
Since  $\varphi_i$ is pseudo-Anosov for  $i \leq j$, 
${\rm ent} \, \varphi_i = || \varphi_i ||_{{\rm T}}$  holds for  $ i \leq j$     
and we have the following corollary,  
\begin{corollary}\label{Cor:Entropy=T} 
The identity 
\begin{equation*}
	{\rm ent} \, \varphi = || \varphi ||_{{\rm T}}
\end{equation*}
holds for any  $\varphi \in \mathcal{M}$. 
\end{corollary} 
\noindent
This should be well known to the experts, 
however it seems to have not been in the literature so far.


\section{Proofs} 

\subsection{Setting up}\label{Subsection:3.1} 

If  $j = 0$, 
then the entropy of  $\varphi$, 
the Weil-Petersson translation distance of  $\varphi$  and the Gromov norm  of  $N_{\varphi}$  
all are obviously zero and we are done for both theorems.  
Thus assume that  $j \geq 1$.  
Then by the theorem of Brock in  \cite{Brock},   
we have a set of inequalities, 
\begin{equation}\label{Eq:BrockInequality}
	C_i^{-1} || \varphi_i ||_{{\rm WP}} \leq {\rm vol} \, N_{\varphi_i} 
		\leq C_i || \varphi_i ||_{{\rm WP}}  \qquad i = 1, 2, \cdots, j. 
\end{equation} 
Each  $C_i$  depends only on the topology of the surface  $S_i$.  
Now, 
there are only finitely many topologies which can appear for 
surfaces obtained by splitting  $\varSigma$  along a system of 
essential curves.  
Thus, 
there is a constant  $F_1= F_1(g, n)$  which bounds the Brock constant 
for a surface appearing in any essential splitting of  $\varSigma$.  
In particular,  
$F_1$  depends only on  $g, n$  and does not on any particular splitting of  $\varSigma$.  
Replacing  $C_i$  by  $F_1$  in (\ref{Eq:BrockInequality}).    
summing all inequalities in  (\ref{Eq:BrockInequality}),   
replacing the middle term by the identity  (\ref{Eq:Gromov=Vol}), 
and we obtain 
\begin{equation}\label{Eq:Comparison1}
	F_1^{-1} \, \sum_{i = 1}^j || \varphi_i ||_{{\rm WP}} 
	\leq v_3 \, || N_{\varphi} ||_{{\rm Gr}} 
	\leq F_1 \sum_{i=1}^j ||\varphi_i ||_{{\rm WP}}.  
\end{equation}

\subsection{Proof of Theorem \ref{Thm:1}} 

By Lemma \ref{Lem:WPleqT} and lemma \ref{Lem:TleqWP}, 
we have a set of inequalities, 
\begin{equation*}\label{Eq:Comparison2}
	E_i \, {\rm ent} \, \varphi_i 
	= E_i \, || \varphi_i ||_{{\rm T}} 
	\leq  || \varphi_i ||_{{\rm WP}} 
	\leq D_i \, || \varphi_i ||_{{\rm T}} 
	= D_i \, {\rm ent} \, \varphi_i 
	\qquad i = 1, 2, \cdots, j.     
\end{equation*}  
By the same reasoning as in subsection \ref{Subsection:3.1}, 
we can find constants  $F_2 = F_2(g, n, \varepsilon)$  
and  $F_3 = F_3(g, n)$  which do not depend on any splitting of  $\varSigma$  
such that the equalities,  
\begin{equation}\label{Eq:Comparison2}
	F_2 \sum_{i=1}^j {\rm ent} \, \varphi_i 
	\leq \sum_{i = 1}^j || \varphi_i ||_{{\rm WP}} 
	\leq F_3 \sum_{i=1}^j {\rm ent} \, \varphi_i,   
\end{equation}  
hold.  
The left inequality is established only under the bounded geometry 
condition controlled by  $\varepsilon$.  

On the other hand, 
since  $j$,  which is the number of pseudo-Anosov components,  
is at most  $n+ 2g-2$, 
there is an obvious constant  $F_4 = F_4(g, n) > 0$  which depends  only on  $g, n$  
such that the inequalties, 
\begin{equation*}\label{Eq:Comparison4}
	F_4^{-1} \max_{1 \leq i \leq j} {\rm ent} \, \varphi_i 
	\leq \sum_{i = 1}^j {\rm ent} \, \varphi_i 
	\leq 	F_4 \max_{1 \leq i \leq j} {\rm ent} \, \varphi_i,   
\end{equation*}  
hold.     
Thus, 
by letting  $F_5 = F_2 F_4^{-1}$  and  $F_6 = F_3 F_4$,  
we have 
\begin{equation}\label{Eq:Comparison5}
	F_5 \, {\rm ent} \, \varphi 
	\leq \sum_{i = 1}^j || \varphi_i ||_{{\rm WP}} 
	\leq F_6 \, {\rm ent} \, \varphi.     
\end{equation}  
Theorem \ref{Thm:1} is immediate from 
the comparisons (\ref{Eq:Comparison1})  
and  (\ref{Eq:Comparison5}).

\subsection{Proof of Theorem \ref{Thm:2}} 
Again since  $j \leq n+2g-2$,  
there is an obvious constant  $F_7 = F_7(g,n)$  which does not 
depend on any splitting of  $\varSigma$  such that  
\begin{equation}
	F_7^{-1} \sqrt{\sum_{i=1}^j || \varphi_i ||_{{\rm WP}}^2}   
	\leq\sum_{i=1}^j ||\varphi_i ||_{{\rm WP}} 
	\leq F_7 \sqrt{\sum_{i=1}^j || \varphi_i ||_{{\rm WP}}^2}.   
\end{equation} 
This together with the identity  ({\ref{Eq:WPnorm})
and the comparison (\ref{Eq:Comparison1})  
immediately imply Theorem \ref{Thm:2}. 
\qed

%
%

\end{document}